\documentclass[journal, a4paper, onecolumn, 12pt]{IEEEtran}

\ifCLASSINFOpdf
\else
\fi
\usepackage{amsmath,amssymb,amsthm,latexsym, amscd,epsfig,epstopdf}
\usepackage{epsfig,graphicx,graphics,subfig,pdfpages, setspace, wrapfig}

\usepackage{amsthm}
\usepackage[utf8x]{inputenc}
\usepackage{amssymb,pdfsync,epstopdf,verbatim,gensymb}
\usepackage{nccmath}
\usepackage{textcomp}
\usepackage{amsmath}
\usepackage{float}
\usepackage{subfig}
\usepackage{amsfonts}
\usepackage{epsfig,color,pbox}

\newtheorem{theorem}{Theorem}[section]

\newcommand{\I}{i}
\newcommand{\D}{d}

\newcommand{\wh}{\widehat}

\newcommand{\lb}{\left(}

\newcommand{\rb}{\right)}

\newcommand{\Oc}{\mathcal{O}}

\newcommand{\Beq}{\begin{equation}}
\newcommand{\Eeq}{\end{equation}}
\newcommand{\beq}{\begin{equation*}}
\newcommand{\eeq}{\end{equation*}}
\newcommand{\bal}{\begin{align}}
\newcommand{\eal}{\end{align}}

\newcommand{\bp}{\begin{prob}}
\newcommand{\ep}{\end{prob}}
\newcommand{\bpr}{\begin{proof}}
\newcommand{\epr}{\end{proof}}

\newcommand{\bel}[1]{\begin{equation}\label{#1}}
\newcommand{\ee}{\end{equation}}

\graphicspath{{./figures/}}

\begin{document}
\title{Numerical inversion of a broken ray transform arising in single scattering optical tomography}
\author{{Gaik Ambartsoumian and Souvik Roy}
\thanks{G. Ambartsoumian was supported in part by US NSF Grant DMS 1109417 and Simons Foundation Grant 360357. S. Roy was supported by the AIRBUS Group Corporate Foundation Chair in Mathematics  of Complex Systems established in TIFR/ICTS, Bangalore, and the University of Texas at Arlington.}
\thanks{G. Ambartsoumian is an Associate Professor at the Department of Mathematics, University of Texas at Arlington, 76019, USA (email: gambarts@uta.edu)}
\thanks{S. Roy was a Visiting Scholar at the Department of Mathematics, University of Texas at Arlington, 76019, USA while carrying out this work. He is currently a Postdoctoral Researcher at Department of Mathematics, International Centre for Theoretical Sciences, 560089, India (email: souvik.roy@icts.res.in) }}
\markboth{}%
{Numerical Inversion of Broken Ray}
\maketitle

\begin{abstract}
 The article presents an efficient image reconstruction algorithm for single scattering optical tomography (SSOT) in circular geometry of data acquisition. This novel medical imaging modality uses photons of light that scatter once in the body to recover its interior features. The mathematical model of SSOT is based on the broken ray (or V-line Radon) transform (BRT), which puts into correspondence to an image function its integrals along V-shaped piecewise linear trajectories. The process of image reconstruction in SSOT requires inversion of that transform.
 We implement numerical inversion of a broken ray transform in a disc with partial radial data. Our method is based on a relation between the Fourier coefficients of the image function and those of its BRT recently discovered by Ambartsoumian and Moon. The numerical algorithm requires solution of ill-conditioned matrix problems, which is accomplished using a half-rank truncated singular value decomposition method. Several numerical computations validating the inversion formula are presented, which demonstrate the accuracy, speed and robustness of our method in the case of both noise-free and noisy data.\end{abstract}



\IEEEpeerreviewmaketitle

\section{Introduction}

Optical tomography uses measurements of light that propagates and scatters inside a body to recover its interior features. If the body part under investigation is optically thick, then the photons scatter multiple times before they are registered outside the body. As a result, the standard mathematical model used in such imaging setups is the diffusion approximation of the radiative transport equation. The latter is severely ill-posed, which leads to subpar quality of reconstructed images. In the optically thin bodies (e.g. biological samples in optical microscopy), most of the photons fly through the object without any scattering, and one can use models based on the regular Radon transform for image reconstruction. While the quality of generated images in this setup is high, the thickness condition restricts its applicability to essentially transparent objects. The middle ground between the two cases described above is the imaging of objects with intermediate optical thickness.  In this case, the input optical photons scatter at most once inside the body and their intensity is measured after they leave it.

The mathematical model for image reconstruction in this single scattering optical tomography (SSOT) is based on the inversion of an integral transform, which puts into correspondence to the image function its integrals along V-shaped piecewise linear flight trajectories of scattered photons. This generalized Radon transform is often called a broken-ray transform (BRT) or a V-line Radon transform (VRT).

SSOT was introduced in a series of influential articles \cite{FMS-PhysRev-10, Florescu-Markel-Schotland, FMS-PhysRev-09} by L. Florescu, J. C. Schotland, and V. Markel, where the authors considered that imaging modality in a rectangular slab geometry. We refer the reader to those papers for details about the underlying physics of SSOT, its advantages over traditional optical imaging modalities, and careful reduction of the radiative transport equation to the integral geometric problem of inverting the BRT in the case of predominantly single scattering of photons.

In this paper we consider SSOT in a circular geometry of data acquisition, where the 2D image function is supported in a disc (3D imaging can be accomplished by vertical stacking of 2D slices as in conventional tomography). Similar to the approach introduced in \cite{FMS-PhysRev-10, Florescu-Markel-Schotland, FMS-PhysRev-09} for slab geometry, here we consider the photons entering the image domain normal to its boundary, travelling a certain distance towards the center of the disc, and then scattering under a certain fixed angle. Notice, that we do not assume that the scattering always happens under the same angle, but rather collect only the data that corresponds to such scattering (e.g. using collimated detectors). As a result the BRT is measured along a two-parameter family of broken rays, where one parameter defines the location of the light source, and the second one the distance from the center of the disc to the scattering location (see Fig. \ref{Sketch_1}). Thus, to recover the 2D image function in the disc one needs to invert the BRT that depends on two variables in circular geometry.

The broken ray transform and its generalization to higher dimensions (called conical Radon transforms) are a fairly recent topic of interest in integral geometry. Their significance grew only a few years ago due to their connections to SSOT and Compton scattering imaging. And while a handful of articles have appeared in literature studying the properties and inversion of BRT in various data acquisition geometries (e.g. see \cite{Ambartsoumian-VLine, Ambartsoumian_Moon_broken_ray_article, Gouia_Amb_V-line, Haltmeier-cones, Jung_Moon, Kats_Krylov-13, Krylov_Kats-15, Nguyen_Truong_Grangeat, Truong_Nguyen_V-line}), many theoretical and practical questions still remain unanswered. In particular, in the case of the circular setup of data acquisition and fixed scattering angle, only two results on BRT are known at this point. We discuss both of them in detail below.

In \cite{Ambartsoumian-VLine} it was shown that if the support of the image function is sufficiently away from the circle of ``source-detector'' positions (see the shaded area in Fig. \ref{Sketch_1}), and BRT data is known for all $\beta\in[0,2\pi]$ and $t\in[-R,R]$, then the problem of inverting the BRT can be reduced to the problem of inverting a regular Radon transform with straight lines passing through the support of image function. Here the negative values of $t$ refer to the broken rays that travel the distance $R+|t|$ from the source to the scattering point (i.e. pass through the center of the disc) and then change direction. Since the inversion of the regular Radon transform is extremely well studied and has many efficient numerical implementations, the approach described in \cite{Ambartsoumian-VLine} works well under the assumptions described above. However, those assumptions are very limiting, especially the requirement of knowing BRT for all $t\in[-R,R]$. As it was mentioned before, the distances that the photons travel in the body induce the number of times that they scatter. Hence, it would be a substantial improvement if in the setup above one could recover the image function just using BRT with $t\in[0,R]$.

In \cite{Ambartsoumian_Moon_broken_ray_article} it was shown that the image function $f(\rho,\phi)$ can indeed be recovered from its BRT $g(\beta,t)$ in circular geometry using just $t\in[0,R]$, $\beta\in[0,2\pi]$ and without the additional restriction on the support of $f$ required in \cite{Ambartsoumian-VLine}. The tradeoff is the complexity of the inversion process in this case. Albeit the analytical inversion formula derived in \cite{Ambartsoumian_Moon_broken_ray_article} is exact, it is based on recovering the Fourier coefficients $f_n(\rho)$ of $f$ through the Fourier coefficients $g_n(t)$ of $g$ using the Mellin transform and its inversion. The numerical implementation of that formula is extremely complicated and was not done in  \cite{Ambartsoumian_Moon_broken_ray_article}.

In this paper, we numerically invert the BRT in the setup described above using the relations between Fourier coefficients of $f$ and $g$ discovered in \cite{Ambartsoumian_Moon_broken_ray_article}. These relations are in the form of integral equations, where the kernel depends on a ratio of the two variables. We solve those by adopting a numerical method given in  \cite{VPVS} and combining it with a truncated singular value decomposition to recover the Fourier coefficients of $f$ from the BRT data $g$.

The rest of the article is organized as follows. Section \ref{Theory_background} gives the relevant theoretical background recalling the integral equations relating $f_n(\rho)$ and $g_n(t)$ based on which the numerical simulations in this paper are performed. Section \ref{Num-algo-section} describes the numerical algorithm for solving the integral equations. In Section \ref{Num_Simulations}, we present the results of the numerical simulations. Section \ref{Summary} concludes the paper with a summary and final remarks.

\section{Theoretical background}\label{Theory_background} 

\begin{figure}
{\centering
    \includegraphics[width=0.5\textwidth]{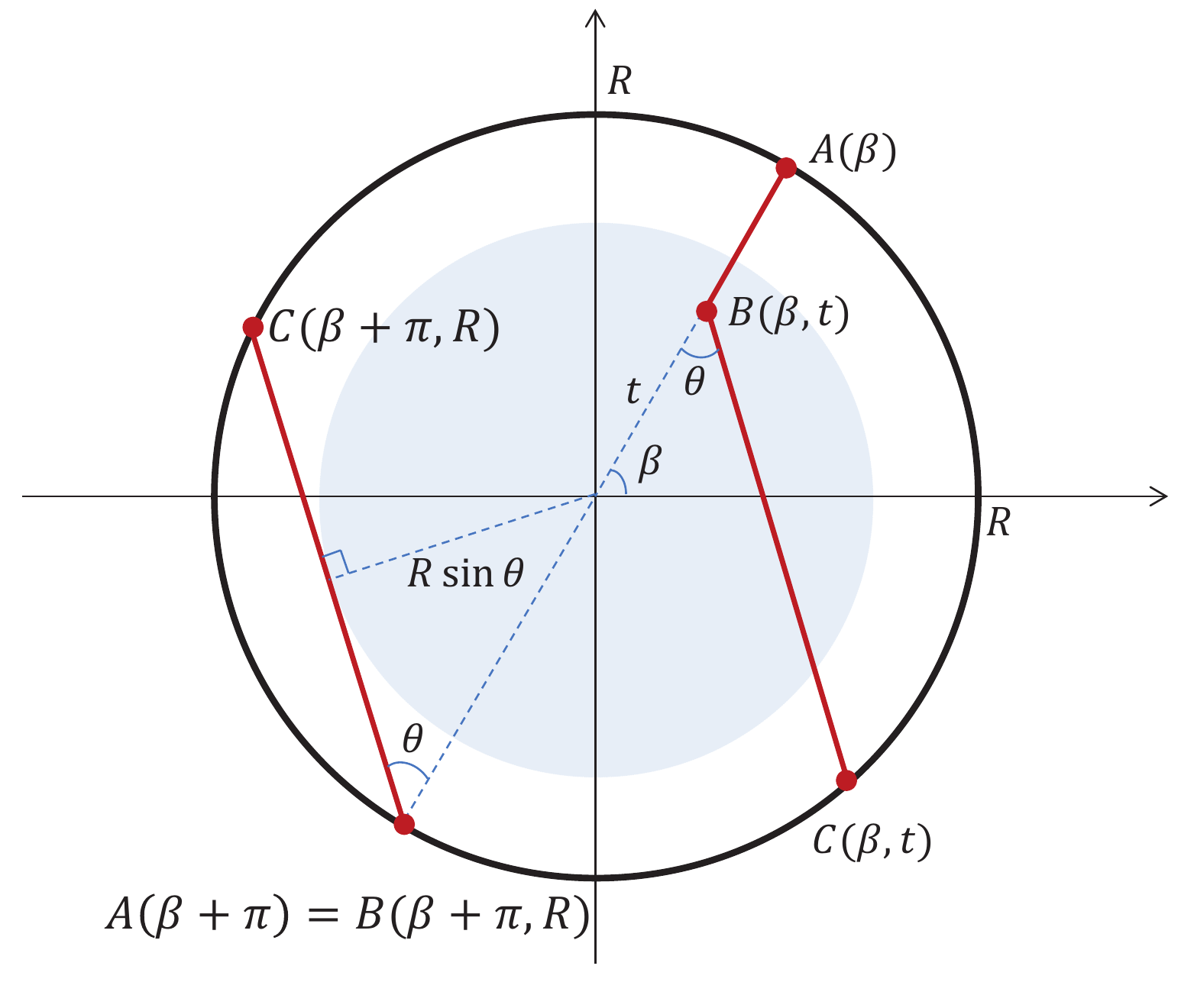}
    \caption{A sketch of BRT in a disc and the notations. Each broken ray is parametrized by two variables: $\beta\in[0,2\pi]$ describing the light source location $A(\beta)$, and $t\in[0,R]$} describing the scattering point $B(\beta, t)$. The detector measuring the intensity of scattered light is located at $C(\beta,t)$.
    \label{Sketch_1}}
\end{figure}

Let the function $f(\rho,\phi)$ be defined inside a disc $D(0,R)$ of radius $R$ centered at the origin, and let $\theta\in (0,\pi/2)$ be a fixed angle. Let $BR(\beta,t)$ denote the broken ray that emits from the point $A(\beta)=(R\cos\beta,R\sin\beta)$ on the boundary of $D(0,R)$, travels the distance $d=R-t$ along the diameter to the point $B(\beta,t)$, then breaks into another ray under the obtuse angle $\pi-\theta$ arriving at point $C(\beta,t)$ (see Fig. \ref{Sketch_1}).

   The broken ray transform of the function $f$ is defined as the integral
\begin{equation}\label{broken_ray_definition}
   \mathcal{R}f(\beta,t) = \int_{BR(\beta,t)} f~ds,~\beta\in[0,2\pi],t\in[0,R],
\end{equation}
of $f(\rho,\phi)$ along the broken ray $BR(\beta,t)$ with respect to linear measure $ds$.

The transform $\mathcal{R}f(\beta,t)$ with radially partial data (i.e. $t\in[0,R]$ instead of $t\in[-R,R]$) was first considered in \cite{Ambartsoumian_Moon_broken_ray_article}. The authors of that paper presented an explicit inversion formula, which was based on the relation between
the Fourier coefficients of $f$ and $g=\mathcal{R}f$, which we recall below.

Expanding $f(\rho,\phi)$ and $g(t,\beta)$ into a Fourier series, we obtain the following
\[
f(\rho,\phi)=\sum\limits_{n=-\infty}^{\infty} f_{n}(\rho)e^{i n\phi},\quad g(t,\beta)=\sum\limits_{n=-\infty}^{\infty} g_{n}(t)e^{\I n\beta}.
\]
We define
\begin{equation}\label{k_n_def}
\begin{aligned}
K_n^1\left({\hat{t}}\right) &= - e^{in \bar{\psi}\left({\hat{t}}\right)}\times 
\frac{1-\hat{t}\cos\left[{\bar{\psi}\left({\hat{t}}\right)}\right]+\hat{t}^2\sin\left[{\bar{\psi}\left({\hat{t}}\right)}\right]\frac{\sin\theta}{\sqrt{1-\hat{t}^2\sin^2\theta}}}{\sqrt{1+\hat{t}^2-2\hat{t}\cos\left[{\bar{\psi}\left({\hat{t}}\right)}\right]}}
\end{aligned}
\end{equation}
and
\begin{equation}\label{k_n_def_2}
\begin{aligned}
K_n^2\left({\hat{t}}\right) &= (-1)^n e^{in\psi\left({\hat{t}}\right)}\times
\frac{1+\hat{t}\cos\left[{\psi\left({\hat{t}}\right)}\right]+\hat{t}^2\sin\left[{\psi\left({\hat{t}}\right)}\right]\frac{\sin\theta}{\sqrt{1-\hat{t}^2\sin^2\theta}}}{\sqrt{1+\hat{t}^2+2\hat{t}\cos\left[{\psi\left({\hat{t}}\right)}\right]}},
\end{aligned}
\end{equation}
where
$$
\bar{\psi}\left({\hat{t}}\right)=2\theta-\psi\left({\hat{t}}\right) \mbox{ and } \psi\left({\hat{t}}\right)= \arcsin\left({\hat{t}\sin\theta}\right)+\theta.
$$

The relation between the $n^{\text{\tiny th}}$ Fourier coefficient of the function $f$ and the $n^{\text{\tiny th}}$ Fourier coefficient of the broken ray transform $g_{n}$ is then given below by the integral equation

\begin{equation}\label{integral_equation}
\begin{aligned}
g_n(t) = & \int_t^R f_n(\rho)d\rho+\int_{t\sin\theta}^t f_n(\rho) K_n^1\left({\frac{t}{\rho}}\right)d\rho + \int_{t\sin\theta}^R f_n(\rho)K_n^2\left({\frac{t}{\rho}}\right)d\rho.
\end{aligned}
\end{equation}
The absolute values of kernels $K^1_n$ and $K^2_n$ have simple geometric meanings. To explain that, let us split the longer branch $BC$
of the broken ray $BR(\beta,t)$ to two parts: $L_1$ from the scattering point $B(\beta,t)$ to the point of the broken ray
closest to the origin, and $L_2$ from the point closest to the origin to the point $C(\beta,t)$. Then the second multiplier
(written as a fraction) in the definition of $K^i_n$ represents the ratio $ds/d\rho$ between the elementary increments of
distance $ds$ along $L_i$ and the corresponding increment of the polar radius $d\rho$.

We now need to solve the integral equation (\ref{integral_equation}). Using the fact that
$$f_n(\rho) = 0,\quad \forall \rho> R,$$
 (\ref{integral_equation}) can be rewritten as
\begin{equation}\label{integral_equation2}
\begin{aligned}
g_n(t) = & \int_t^R f_n(\rho)d\rho+\int_{t\sin\theta}^t f_n(\rho) K_n^1\left({\frac{t}{\rho}}\right)d\rho\\
 +& \int_{t\sin\theta}^R f_n(\rho)K_n^2\left({\frac{t}{\rho}}\right)d\rho\\
=&\int_{t\sin\theta}^R f_n(\rho)K_n\left(\frac{t}{\rho}\right)\D \rho,
\end{aligned}
\end{equation}
where
\begin{equation}
K_n\left({\hat{t}}\right)=\left\{
  \begin{array}{l l}
     K_n^1\left({\hat{t}}\right)+K_n^2\left({\hat{t}}\right),& \quad \text{if }1< \hat{t}< 1/\sin\theta,\\
    1+K_n^2\left({\hat{t}}\right),& \quad \text{if } 0 \leq \hat{t} \leq 1\\
  \end{array} \right.
\end{equation}
and $K_n^1\left({\hat{t}}\right),K_n^2\left({\hat{t}}\right)$ are given by (\ref{k_n_def}) and (\ref{k_n_def_2}).

Notice, that $\hat{t}_1=1/\sin(\theta)$ corresponds to the point of the broken ray that is the closest to the origin,
while $\hat{t}_2=1$ corresponds to the the scattering point $B$ on $L_1$ and its symmetric point on $L_2$.
Formulas (\ref{k_n_def}) and (\ref{k_n_def_2}), as well as the geometric interpretation given above, show that
the function $K_n(\hat{t})$ is infinitely smooth on the interval $[0,\hat{t}_1)$, except at point $\hat{t}_2$, where it has
a jump discontinuity of size $1/\cos(\theta)-1$. One can also notice that, $K_n(\hat{t})$ blows up to infinity as $\hat{t}$ approaches $\hat{t}_1$
(due to $\sqrt{\hat{t}_1-\hat{t}}$ in the denominator), but that singularity is integrable.

Let us rewrite equation (\ref{integral_equation2}) as follows
\begin{equation}\label{final_integral_eq}
g_n(t) = \int_{t\sin\theta}^R f_n(\rho)\;K_n\left(\frac{t}{\rho}\right)\D \rho.
\end{equation}

A simple change of variables $t=e^{\tau}$ and $\rho=e^r$ will transform it to a convolution type integral equation
\begin{equation}\label{convolution_equation}
G_n(\tau)=\int\limits_{\tau+\ln\sin{\theta}}^{\ln{R}} F_n(r)\, k_n(\tau-r)\, dr,
\end{equation}
where $G_n(\tau)=g_n(e^{\tau})$, $F_n(r)=e^r f_n(e^r)$, and $k_n(x)=K_n(e^x)$.

Assume that $f$ is an infinitely differentiable function compactly supported in
$A(\varepsilon,R)$ (the annulus centered at the origin, with inner radius $\varepsilon>0$ and exterior radius $R$).
Then equation (\ref{final_integral_eq}) can be considered for $t\in[\varepsilon/\sin(\theta), R]$, where $\varepsilon>0$
is an arbitrarily small number. As a result, (\ref{convolution_equation}) becomes a convolution type integral equation, where $F_n$ is infinitely differentiable and compactly supported, and $k_n$ is locally integrable and compactly supported.
Such an equation can be solved by taking a Fourier transform of both sides, dividing by the Fourier transform of the kernel,
and taking an inverse Fourier transform.
As a result, we get the following

\begin{theorem}[Existence and uniqueness of solution]\label{Existence_Uniqueness_Mellin}
Let $f$ be a $C^\infty$ function with support inside the annulus $A(\varepsilon,R)$. Then equation (\ref{final_integral_eq})
has a unique solution $f_{n}(t)\in C^\infty([\varepsilon, R])$.
\end{theorem}

The smoothness of $f_n$ follows from the smoothness of $f$ and uniqueness of the solution of (\ref{final_integral_eq}).
Since $\varepsilon>0$ in the theorem above can be taken arbitrarily small, in our numerical experiments described below
we can assume the hypothesis of the theorem is satisfied simply by using a discretization that avoids the broken ray passing
through the origin.

\section{Numerical Algorithm}\label{Num-algo-section}In this section, we describe the numerical scheme used to solve the integral equation (\ref{final_integral_eq}).

\subsection{Fourier coefficients of the broken ray Radon data in the angular variable}\label{rd}

The function $g(t,\beta)$ is real-valued in the angular variable $\beta$. One could perform the standard FFT on the discrete sequence of values $\lbrace{g(t,\beta_N)}\rbrace_N$. However, there is a computationally more efficient way of performing the FFT on such a real-valued function. Here the real data of length $N$ is split into two equal halves and a complex data of length $N/2$ is created. On this complex data, FFT is performed and then converted back to FFT of real data. Therefore, for large $N$, almost half of the arithmetic operations can be saved by performing the FFT on $N/2$
complex numbers instead of treating the real sequence as consisting of $N$ complex numbers. 

 Thus we compute the modified discrete fast Fourier transform (FFT) of $g(t,\beta)$ in $\beta$ for a fixed $t\in [0,R]$ based on the Cooley-Tukey algorithm (see \cite{MR1201159}) as follows

\begin{enumerate}
 \item
 Let $\{\beta_1, \beta_2,\cdots,\beta_N\}$ be a discretization of $\beta$, where $N$ is even. We break the array $g(t,\beta_{k})$ for $1\leq k\leq N$ into two equal length arrays for odd and even numbered indices. Thus, we define $A = \{g(t,\beta_{2j-1})\}$ and $B = \{ g(t, \beta_{2j})\}$ for $j = 1,2,\cdots,N/2$.
 \item Next we create a complex array $h^{c}_{t}(j) = A(j) + \I B(j),\quad j = 1,2,\cdots,N/2$.
 \item We now perform a discrete FFT on $h^c_t$ to get $\wh{h^{c}_t}(n), \quad n = 1, 2,\cdots, N/2$.

 \item Thus the Fourier series of $g$ in the $\beta$ variable is given as follows
 \begin{equation}
 \begin{aligned}
 g_n(t) =
 \begin{cases}
  \frac{1}{2}\Big{\{}\lb \wh{h^c_t}(n)+ \overline{\wh{h^c_t}(\frac{N}{2}-n+2)}\rb\\
  -\I\lb\wh{h^c_t}(n)- \overline{\wh{h^c_t}(\frac{N}{2}-n+2)}\rb\cdot e^{\frac{2\pi{i}(n-1)}{N}}\Big{\}},\\
  \hspace{30mm}\mbox{for } n = 1,\cdots,\frac{N}{2}+1 \\
  \wh{h^c_t}(N -n+2),~ n = \frac{N}{2}+2,\cdots,N.
  \end{cases}
\end{aligned}
\end{equation}

\end{enumerate}

\subsection{Trapezoidal product integration method }
Using the values of $g_n(t)$ obtained in the previous section, we want to solve the integral equation (\ref{final_integral_eq}). Under the assumptions of Theorem \ref{Existence_Uniqueness_Mellin} and using the method of Mellin transforms, the authors in \cite{Ambartsoumian_Moon_broken_ray_article} provide an analytical inversion formula for the broken ray transform with radially partial data. However, that exact inversion formula is numerically unstable.
One of the reasons for the instability is the fact that the function $K_n(\hat{t})$ is unbounded in the neighborhood of the point $1/\sin(\theta)$. And while its Mellin transform is well-defined as an improper integral, its numerical computation is highly unstable.
Therefore, we approach the numerical inversion problem by solving \eqref{final_integral_eq} directly.
We use the so-called trapezoidal product integration method proposed in  \cite{VPVS,Weiss_Product_Integration_Paper}. We briefly sketch this method below.

\begin{figure}
\centering
    \includegraphics[width=0.5\textwidth]{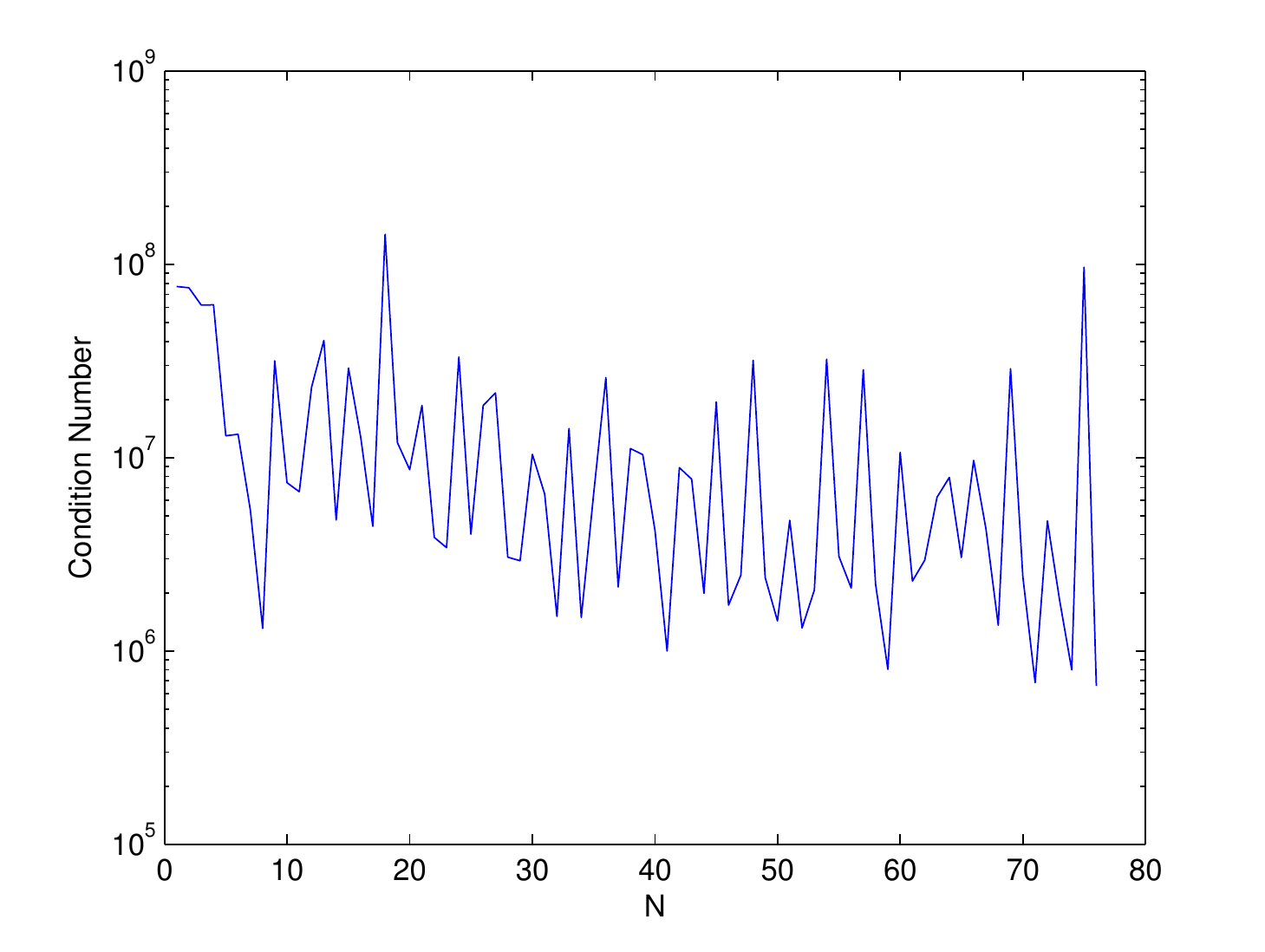}
    \caption{Plot of condition number of $A_n$ for $n\in[0,76]$  }
    \label{fcond}
\end{figure}

Let $M$ be a positive integer and $t_{l}=lh,  l = 0,\hdots,M$ and $ h = R/M$ be a discretization of $[0,R]$.
Choose $i \in \lbrace 0,\hdots, M\rbrace$. Thus from (\ref{final_integral_eq}) we have
\begin{equation}\label{disc_eq}
g_n(t_i) = \int_{t_i\sin\theta}^R f_n(\rho)K_n\left(\frac{t_i}{\rho}\right)\D \rho,
\end{equation}
We choose an index $l$ such that
$$
t_{l} \leq t_i\sin\theta < t_{l+1}.
$$
(If there exists no such $l$ satisfying $t_{l} \leq t_i\sin\theta$, we choose $l=0$).
We then approximate (\ref{disc_eq}) as
\begin{equation}
g_n(t_i) = \int_{t_l}^R f_n(\rho)K_n\left(\frac{t_i}{\rho}\right)\D \rho.
\end{equation}
In the sub-interval $[t_{k},t_{k+1}]$, we approximate $f_{n}(\rho)K_{n}\left(\frac{t_i}{\rho}\right)$ by a linear function taking the values $f_{n}(t_{k})K_{n}\left(\frac{t_i}{t_k}\right)$ and $f_{n}(t_{k+1})K_{n}\left(\frac{t_i}{t_{k+1}}\right)$ at the endpoints $t_{k}$ and $t_{k+1}$, respectively. This is given by
\begin{equation}
\begin{aligned}
f_n(\rho)~K_{n}\left(\frac{t_i}{\rho}\right)\approx& f_{n}(t_{k})K_{n}\left(\frac{t_i}{t_k}\right)\frac{t_{k+1}-\rho}{h}
+f_{n}(t_{k+1})K_{n}\left(\frac{t_i}{t_{k+1}}\right)\frac{\rho-t_k}{h}.
\end{aligned}
\end{equation}
Hence
\begin{equation}
\begin{aligned}
g_{n}(t_{i}) \approx& \sum_{k=l}^{M-1}\int_{t_{k}}^{t_{k+1}}
f_{n}(t_{k})K_{n}\left(\frac{t_i}{t_k}\right)\frac{t_{k+1}-\rho}{h}
+f_{n}(t_{k+1})K_{n}\left(\frac{t_i}{t_{k+1}}\right)\frac{\rho-t_k}{h}\D \rho\\
=& \sum_{k=l}^{M-1}f_{n}(t_{k})K_{n}\left(\frac{t_i}{t_k}\right)\frac{h}{2}+f_{n}(t_{k+1})K_{n}\left(\frac{t_i}{t_{k+1}}\right)\frac{h}{2}.
\end{aligned}
\end{equation}
Thus
\begin{equation}\label{scheme}
\begin{aligned}
\sum_{k=0}^{M-1}a_{i,k}f_{n}(t_{k})=g_n(t_i)
\qquad i = 1,\cdots,M,\\
\end{aligned}
\end{equation}
where
\begin{equation}\label{coefficients}
 a_{i,k} = \left\{
  \begin{array}{l l }
      &\dfrac{h}{2}K_{n}\left(\dfrac{t_i}{t_k}\right),\quad k=l,M,\\
      &h K_{n}\left(\dfrac{t_i}{t_k}\right),\quad l < k < M,\\
      & 0,\quad \mbox{ otherwise.}\\
  \end{array} \right.
\end{equation}
The following theorem states the error estimate for the numerical solution of the integral equation (\ref{final_integral_eq}), which can be proved using the arguments given in \cite[Thm. 7.2]{Linz_book}, once we choose the point of discontinuity $t=1$ as a nodal point.
\begin{theorem}[Error Estimates]\label{Error Estimates}
Let $f_{n}^{\mathrm{exact}}$ be the $C^{3}$ solution of (\ref{final_integral_eq}) in $[0,R]$ and $f_{n}$ be the solution to \eqref{scheme}. Then
\Beq
 \Vert f_{n}^{\mathrm{exact}} - f_{n}\rVert_2= \Oc(h^2),
\Eeq
where $\Vert\cdot\rVert_2$ represents the discrete version of the continuous $L^2$ norm in $[0,R]$ (see for e.g., \cite[Ch. 4]{book_L2norm}).
\end{theorem}
Equation (\ref{scheme}) can be written in matrix form as
\begin{equation}\label{mat}
 A_nF_{n} = \widetilde{g}_{n},
\end{equation}
where  \Beq\label{Discretization of F}F_{n}= \left( \begin{array}{c}
f_{n}(t_{1})  \\
\vdots  \\
f_{n}(t_{M})  \end{array} \right), \qquad
g_{n}= \left( \begin{array}{c}
g_{n}(t_{1})  \\
\vdots  \\
g_{n}(t_{M})  \end{array} \right).
\Eeq
The matrix $A_n$ is defined as
\begin{equation}\label{Matrix equation} A_n(i,k) = a_{i,k},
  \end{equation}
 where $a_{i,k}$ is given by (\ref{coefficients}).
To solve (\ref{scheme}), we need to invert the matrices $A_n$. It turns out that the matrices $A_n$ are ill-conditioned. Fig. \ref{fcond} shows the condition number of $A_n$ for different values of $n$.
It is well known that numerically inverting a matrix with condition number $r$ leads to a loss of $r$ digits of accuracy (see \cite{Hansen_TSVD_Paper}). From Figure \ref{fcond}, we see that the condition numbers of $A_{n}$ is greater than $10^6$ for almost all values of $n$. Thus for the inversion of $A_n$, we use the Truncated Singular Value Decomposition (TSVD) (see \cite{Hansen_TSVD_Paper}) to solve the matrix equation (\ref{mat}).

\begin{figure}
\centering
    \includegraphics[width=0.5\textwidth]{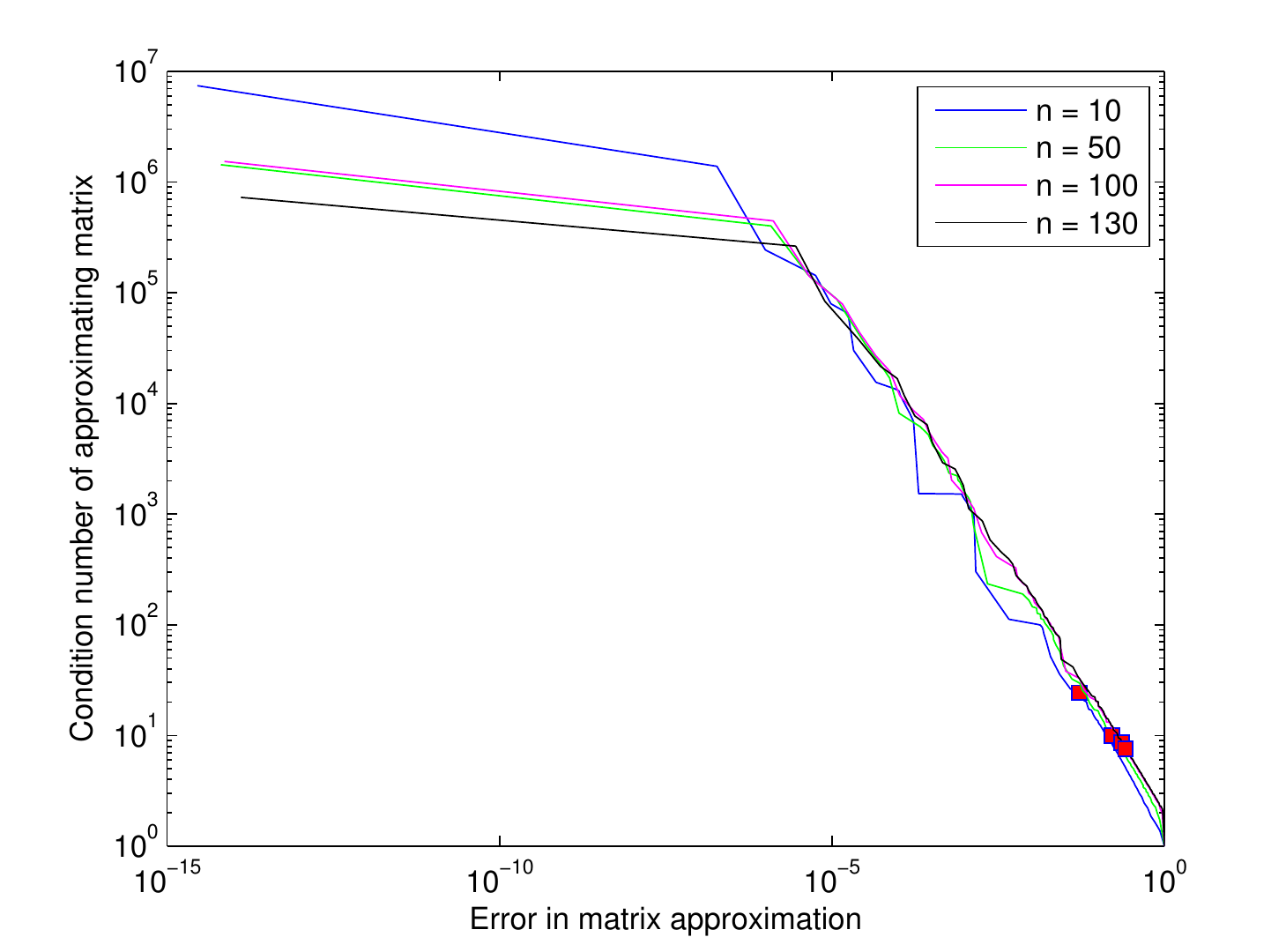}
    \caption{Relation between condition number of $A_{n,r}$ and the error in 2-norm for the original matrix $A_n$, respectively for $n=10,50,100,130.$ The dots in the figure correspond to half-rank approximations.  }
    \label{cond_err}
\end{figure}
\subsection{Truncated singular value decomposition (TSVD)} \label{TSVD}
To solve (\ref{mat}), we first compute the SVD of $A_{n}$. This is given by $A_{n}=UDV^{T}$, where $U$ and $V$ are orthogonal matrices whose columns are the eigenvectors of $A_{n}A_{n}^{T}$ and $A_{n}^{T}A$ respectively. The matrix $D$ is a  diagonal matrix consisting of the singular values of $A_{n}^{T}A_{n}$, which represents the square root of the eigenvalues of $A_{n}^{T}A_{n}$ in descending order represented by $\sigma_r$. We set
 \[ A_{n,r}=UD_{r}V^{T} \quad\mbox{and}\quad A_{n,r}^{-1} = VD^{-1}_{r}U^{T}
 \]
 where $D_{r}$ and $D_{r}^{-1}$ are diagonal matrices with diagonal  entries
 \[
 \begin{array}{ll}
 \lb D_{r}\rb_{ii}=
 \begin{cases}
 &D_{ii} \quad \mbox{if}\quad  i\leq r\\
 &0\quad  \mbox{     otherwise. }
 \end{cases}
 \end{array}
 \]
 \[
 \begin{array}{ll}
 \lb D_{r}^{-1}\rb_{ii}=
 \begin{cases}
 &\frac{1}{D_{ii}} \quad \mbox{if}\quad  i\leq r\\
 &0\quad  \mbox{     otherwise. }
 \end{cases}
 \end{array}
 \]
 Then the matrices $A_{n,r}$ approximates $A_{n}$ for $1\leq r\leq M$. $r$ is the rank of the matrix $A_{n,r}$ (see \cite{Hansen_TSVD_Paper}). We define the matrix 2-norm or the spectral norm of a matrix $A$ of order $n$ as follows
\begin{equation}\label{spec}
\|A\|_2=\max_{|x|_2\neq 0}\frac{|Ax|_2}{|x|_2},
\end{equation}
where $x\in\mathbb{R}^n$ and $|x|_2=(\sum_{i=1}^nx_i^2 )^{1/2}$.

 The condition number of the truncated matrix $A_{n,r}$ is defined to be $\kappa(A_{n,r}) = \dfrac{\sigma_1}{\sigma_r}$. Moreover, the error in approximation of $A$ by $A_{n,r}$ is defined as $\|A-A_{n,r}\|_2$ and is given by $\|A-A_{n,r}\|_2 = \sigma_{r+1}$.

Fig. \ref{cond_err} shows the relation between the condition number of the truncated matrix $A_{n,r}$ and the error $\|A-A_{n,r}\|_2$. A high rank approximation would render the condition number of $A_{n,r}$ to be large, whereas a low rank approximation would lead to loss of information resulting in incomplete reconstruction. In this paper, we have taken half-rank approximations, that is, $r=\frac{M}{2}$.

\section{Numerical Results}\label{Num_Simulations}

We now validate the numerical algorithm proposed in Sec. \ref{Num-algo-section} for the broken ray Radon transforms given in (\ref{broken_ray_definition}). We discretize $\phi \in [0, 2\pi]$ into 150 equally spaced grid points. For discretization of $\rho \in [0, R − \epsilon]$ (we chose $\epsilon = 0.001$), we consider 150 and 400 equally spaced grid points. In all cases we chose $R = 1$. We organize this section as follows. In Sec. \ref{sec:Radon_data}, we describe the procedure of generating the Radon data. In Sec. \ref{sec:test_cases}, we demonstrate our algorithm on two test cases, and analyze the visual features of the reconstructed images. A detailed description of artifacts appearing in the reconstructions is given in Sec. \ref{sec:stability_disc}. Finally, in Sec. \ref{sec:error}, we provide the computational times taken for the simulations and evaluate the relative $L^2$ error percentage between the actual and the reconstructed images.

\subsection{Generating the Radon data}\label{sec:Radon_data}
The first step to validate the numerical algorithm described in Sec. \ref{Num-algo-section} is to generate the Radon data. We divide our domain $D(0,R)$ into $150 \times 150$ pixels. The phantom $f$ to be reconstructed is represented by the information present in the pixels. In the next step, we fix $t$ and $\beta$, as defined in Sec. \ref{Theory_background} and consider the measure of intersection of the broken ray $BR(\beta,t)$ with a pixel. Multiplying the values of the measures obtained with the values in the pixels and then summing up over all pixels gives us the Radon data $Rf(\beta,t)$.

\subsection{Test Cases}\label{sec:test_cases}

We now proceed to the test cases to validate our numerical algorithm. The computations are done in MATLAB 7.14.0.739, on a Intel I7 3.1 GHz. quad core processor with 6 GB RAM.

\subsubsection{ Test Case 1 - Disk containing origin}\label{case_1}
\begin{figure}[H]
  \centering
    \subfloat[]{%
     \label{Disk_origin_phantom_exact} \includegraphics[scale=0.5,keepaspectratio]{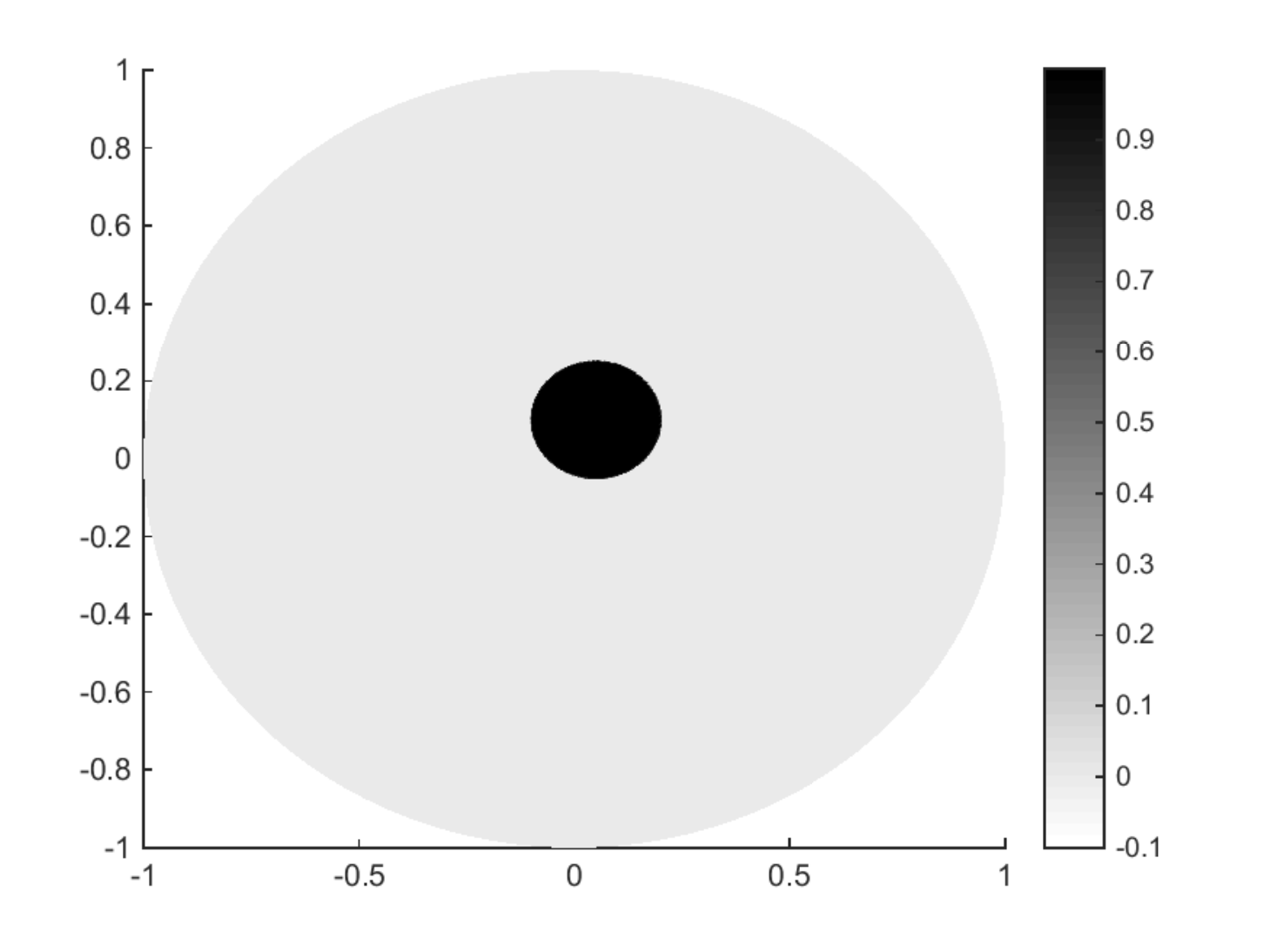}
   }
    \subfloat[]{%
      \label{Reconstruction_disk_origin}\includegraphics[scale=0.5,keepaspectratio]{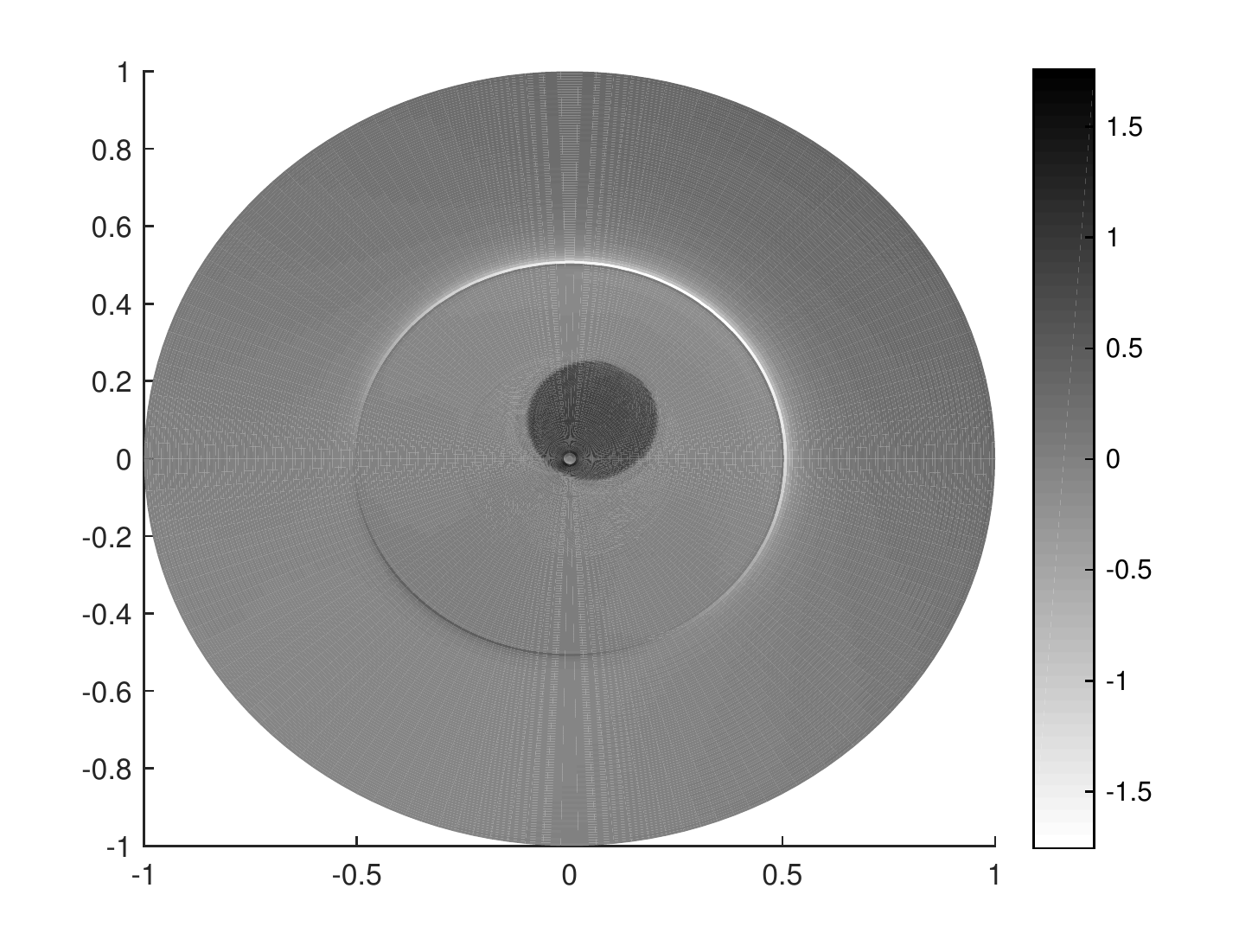}
      }\\
      \subfloat[]{%
      \label{Reconstruction_disk_origin_refined}\includegraphics[scale=0.5,keepaspectratio]{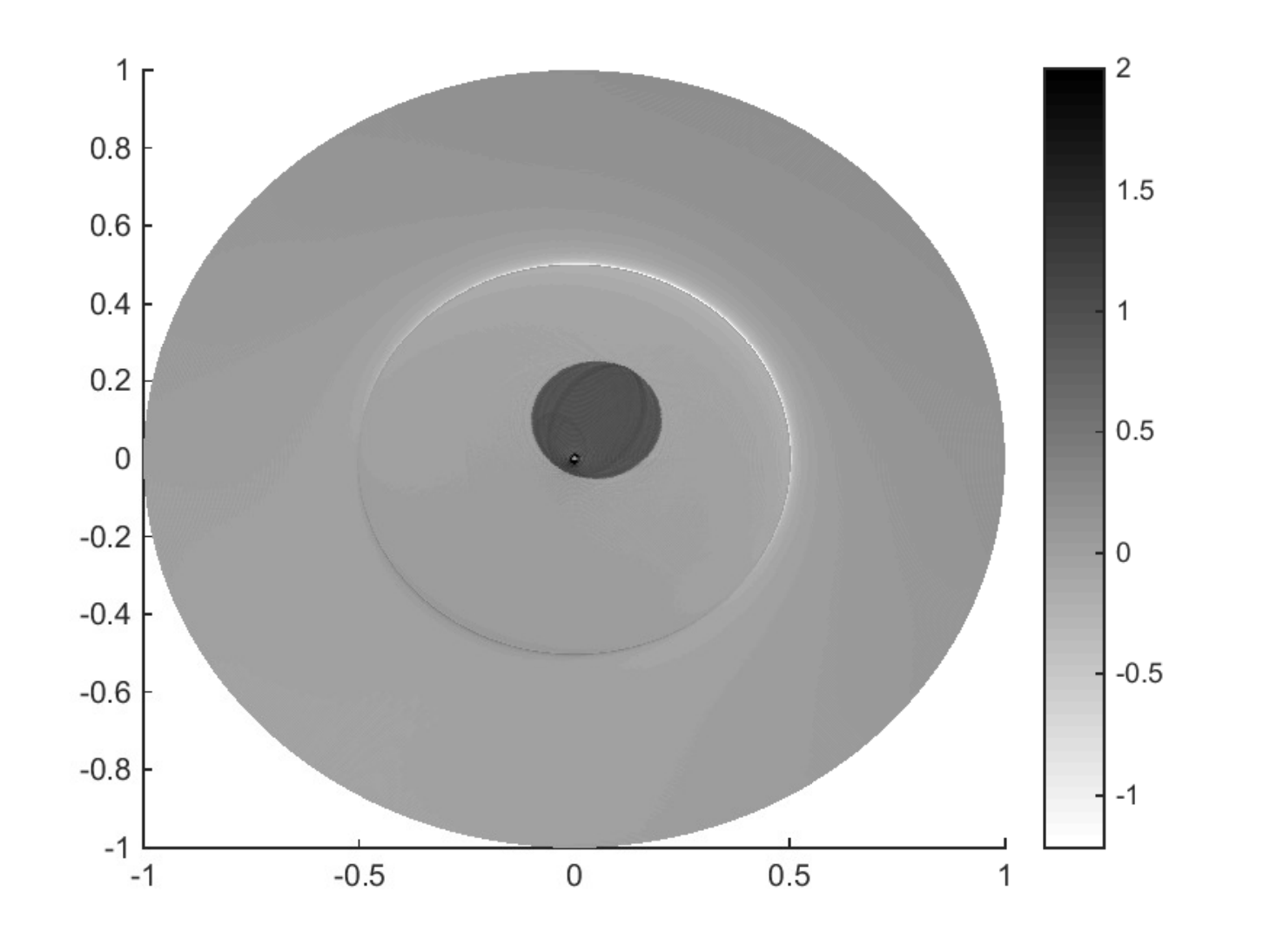}
      }
      \caption{Results for broken ray transform data with breaking angle $\theta=\pi/6$, for a phantom represented by a disk. Figure \ref{Disk_origin_phantom_exact} represents the actual phantom. Figures \ref{Reconstruction_disk_origin}, \ref{Reconstruction_disk_origin_refined} show the reconstructed images with 150 and 400 equally spaced discretizations, respectively.}
      \label{fig:disk_origin}
      \end{figure}

\begin{figure}[H]
  \centering

    \subfloat[]{%
      \label{Reconstruction_disk_origin_pi_4}\includegraphics[scale=0.5,keepaspectratio]{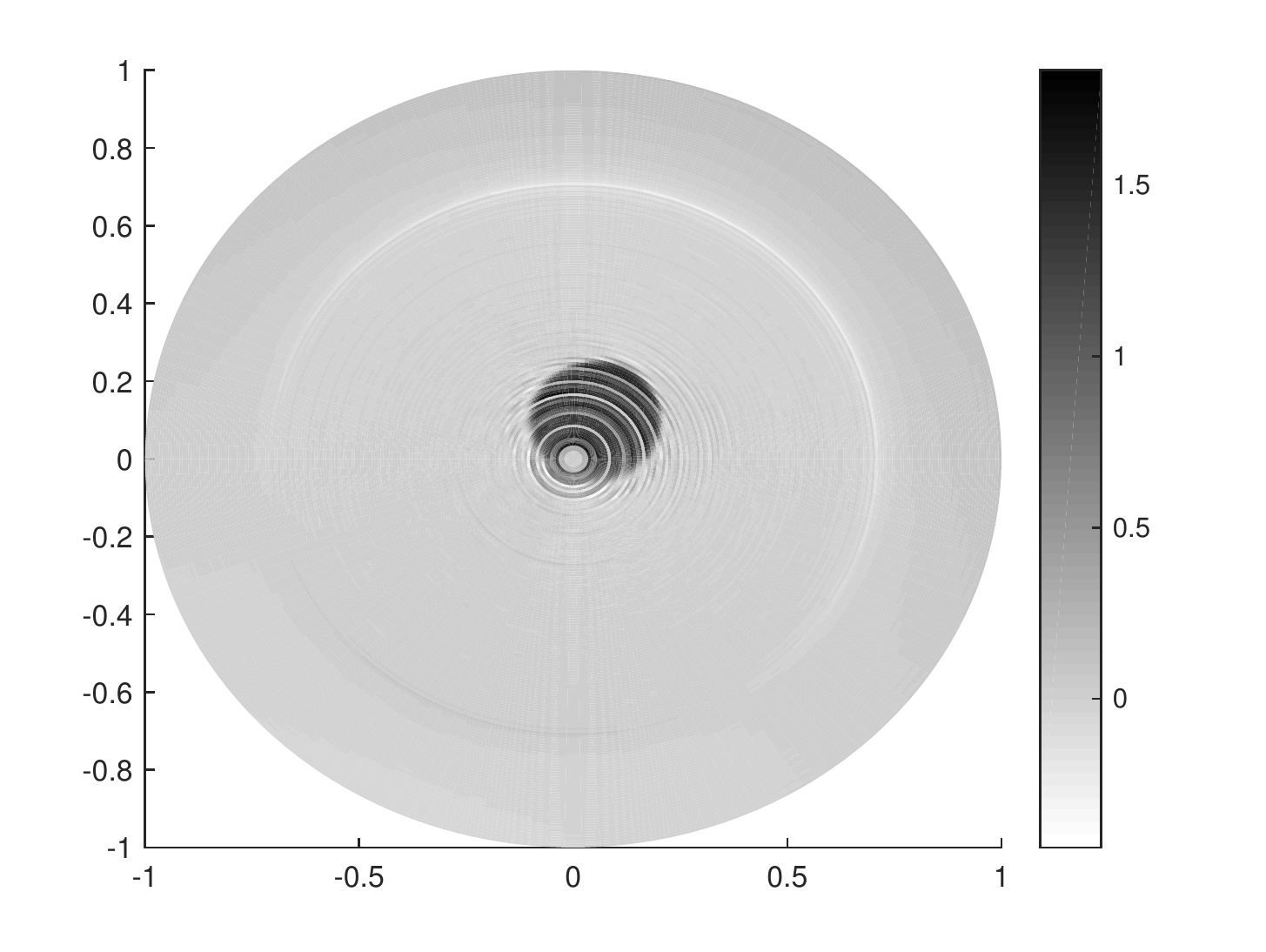}
      }
      \subfloat[]{%
      \label{Reconstruction_disk_origin_pi_4_refined}\includegraphics[scale=0.5,keepaspectratio]{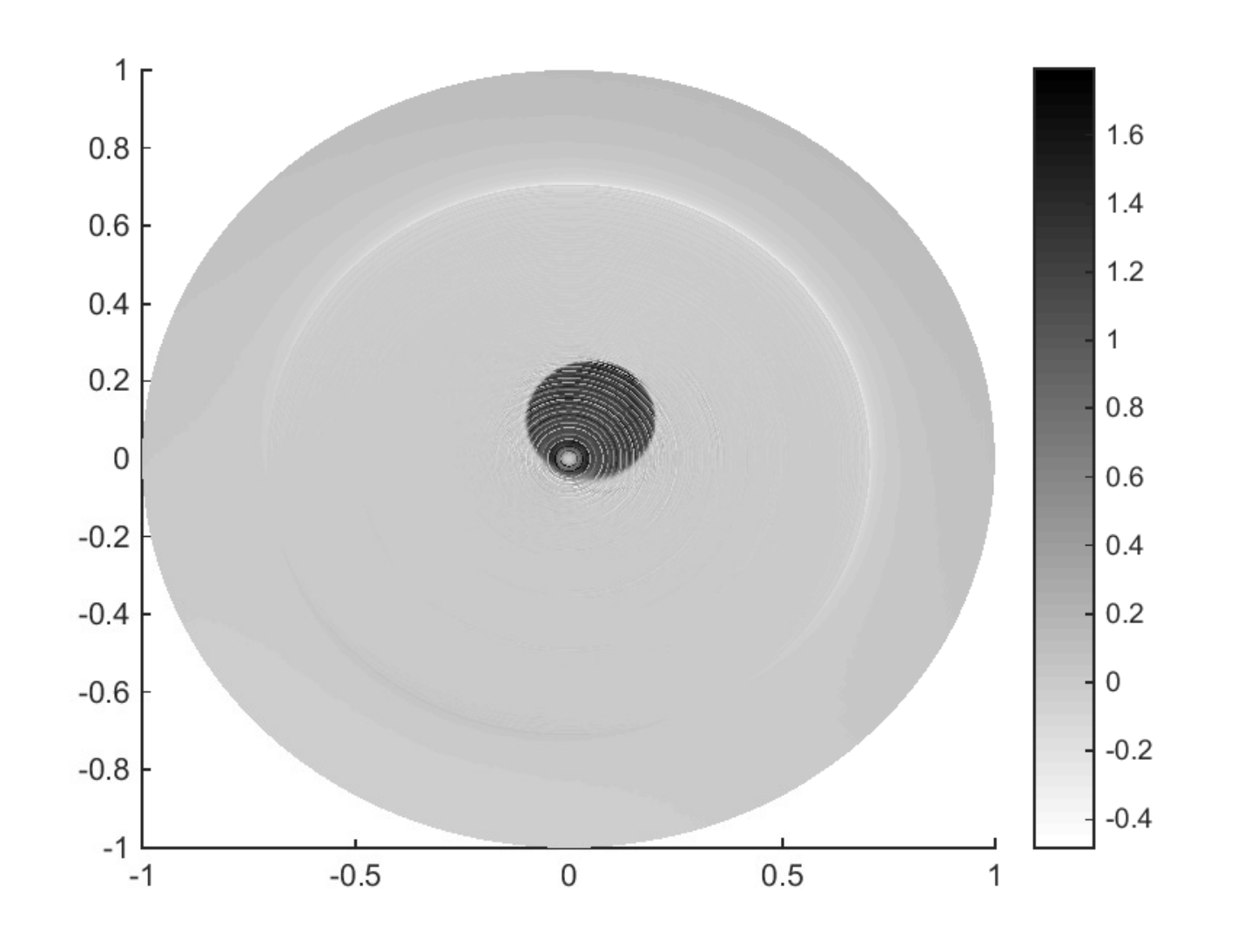}
      }
      \caption{Results for broken ray transform data with breaking angle $\theta=\pi/4$, for a phantom represented by a disk as shown in Figure \ref{Disk_origin_phantom_exact}. Figures \ref{Reconstruction_disk_origin_pi_4}, \ref{Reconstruction_disk_origin_pi_4_refined} show the reconstructed images with 150 and 400 equally spaced discretizations, respectively.}      \label{fig:disk_origin}
      \end{figure}

Fig. \ref{Disk_origin_phantom_exact} shows a phantom represented by a disk centered at $(0.05,0)$ with radius 0.15, thus, containing the origin. Fig. \ref{Reconstruction_disk_origin} and \ref{Reconstruction_disk_origin_refined} show the reconstructions with $\theta=\pi/6$ for 150 and 400 equally spaced discretizations in $\rho$, respectively. Fig. \ref{Reconstruction_disk_origin_pi_4} and \ref{Reconstruction_disk_origin_pi_4_refined} show the reconstructions with $\theta=\pi/4$ for 150 and 400 discretizations in $\rho$ respectively. We see that the phantom is reconstructed fairly well. Not only there is a proper recovery of the shape and the location of the phantom, the values are also well approximated. The visible artifacts along the circle of radius $R\sin\theta$ and at the origin are discussed later in Sec. \ref{sec:stability_disc}.

\begin{figure}
  \centering

    \subfloat[]{%
      \label{Reconstruction_disk_origin_noise}\includegraphics[scale=0.5,keepaspectratio]{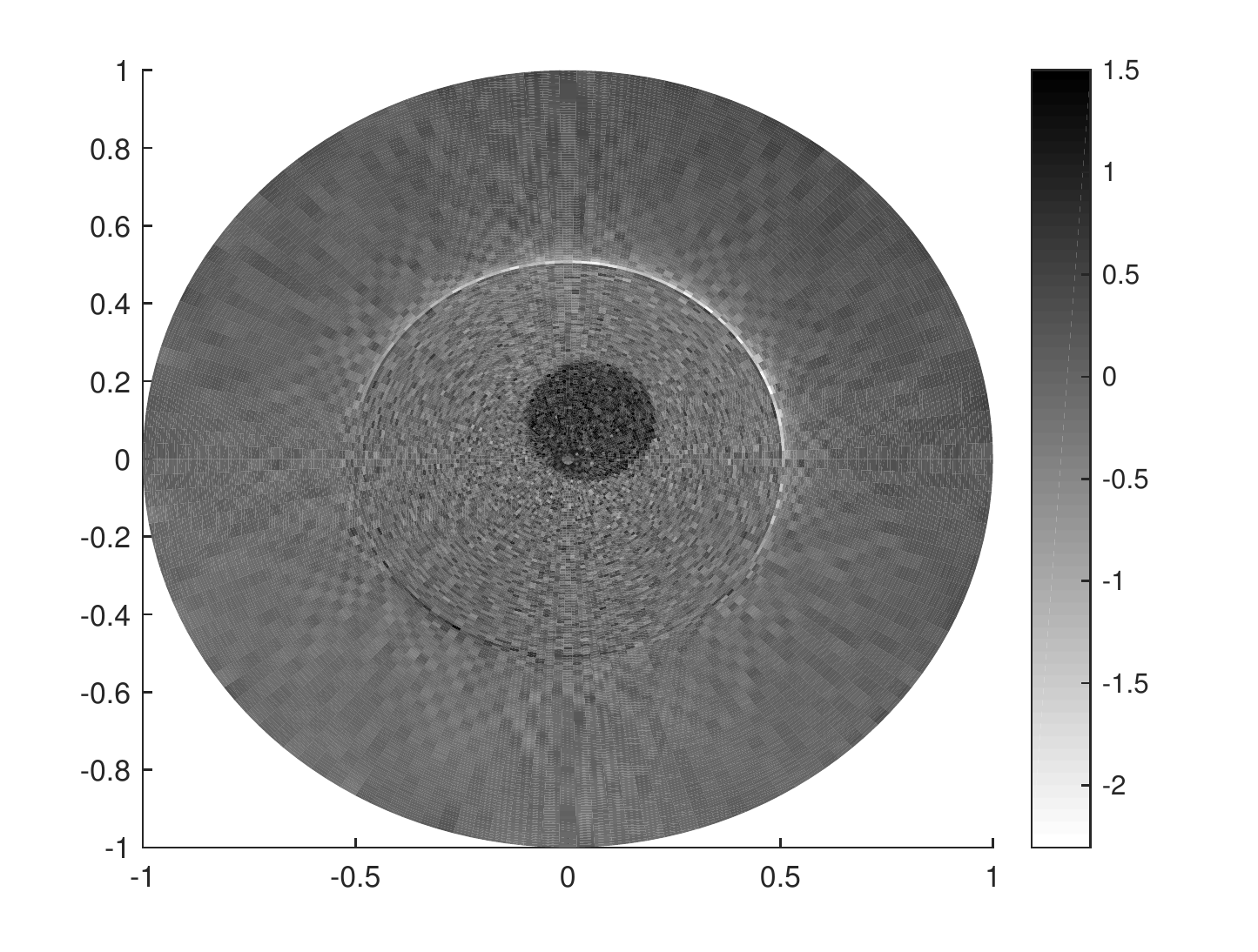}
      }
      \subfloat[]{%
      \label{Reconstruction_disk_origin_refined_noise}\includegraphics[scale=0.5,keepaspectratio]{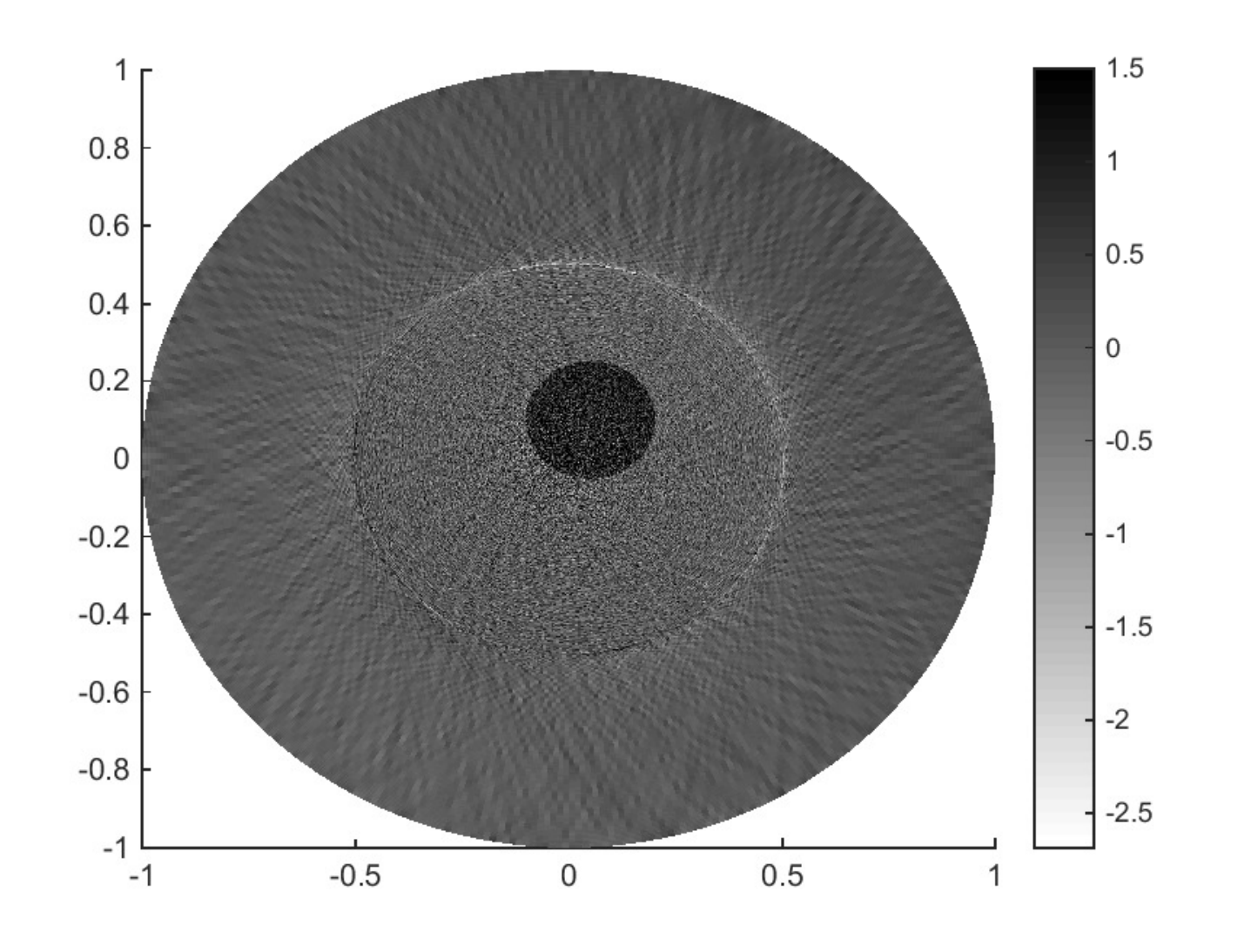}
      }
      \caption{Results for broken ray transform data with breaking angle $\theta=\pi/6$ and with $5\%$ multiplicative Gaussian noise, for a phantom represented by a disk as shown in Figure \ref{Disk_origin_phantom_exact}. Figures \ref{Reconstruction_disk_origin_noise}, \ref{Reconstruction_disk_origin_refined_noise} show the reconstructed images with 150 and 400 equally spaced discretizations, respectively.}
      \label{fig:disk_origin}
      \end{figure}

To demonstrate the robustness of our algorithm, we also tested it for inverting the Radon data with $5\%$ multiplicative Gaussian noise. Fig. \ref{Reconstruction_disk_origin_noise} and \ref{Reconstruction_disk_origin_refined_noise} show the reconstructions  for 150 and 400 equally spaced discretizations in $\rho$, respectively. We again note the good recovery in both cases.

\begin{figure}
  \centering

    \subfloat[]{%
      \label{Less_than_half_rank}\includegraphics[scale=0.5,keepaspectratio]{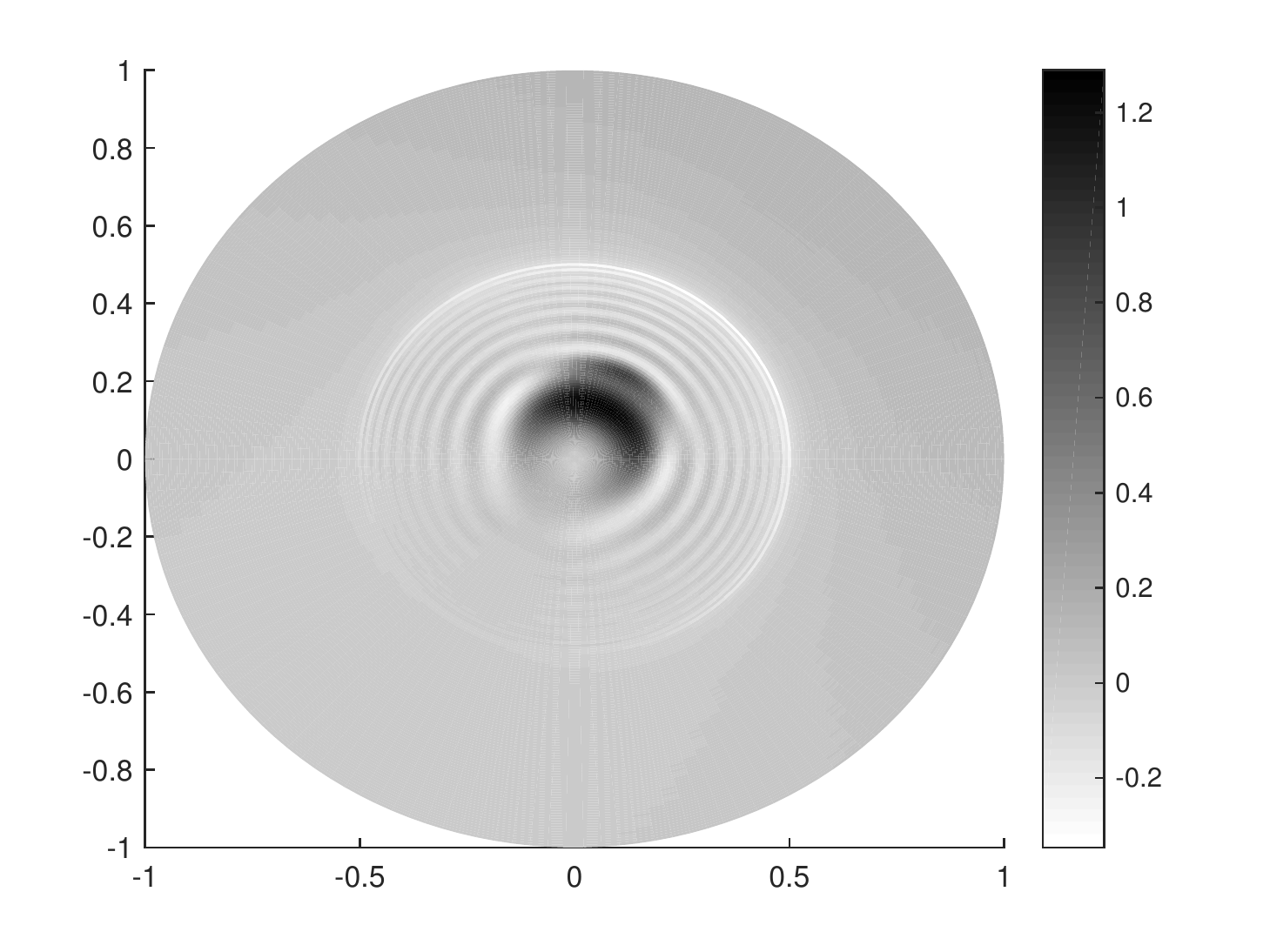}
      }
      \subfloat[]{%
      \label{More_than_half_rank}\includegraphics[scale=0.5,keepaspectratio]{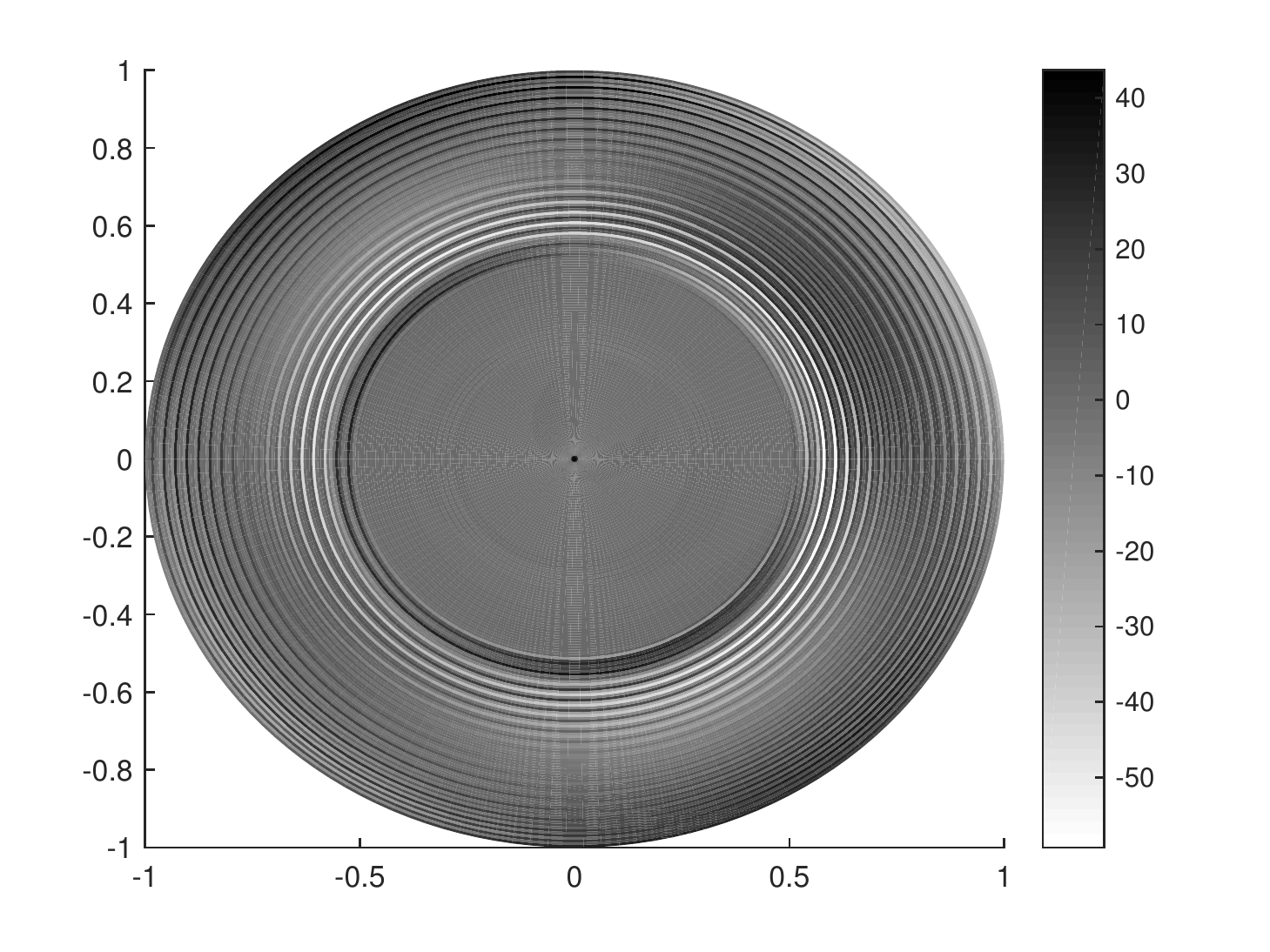}
      }
      \caption{Results for broken ray transform data with breaking angle $\theta=\pi/6$, for a phantom represented by a disk as shown in Figure \ref{Disk_origin_phantom_exact}. Figure \ref{Less_than_half_rank} shows the reconstruction with rank $r=M/8$. Figure \ref{More_than_half_rank} shows the reconstruction with $r=M/1.5$. Figure \ref{Less_than_half_rank} reveals incomplete reconstruction due to loss of data, whereas Figure \ref{More_than_half_rank} reveals blow-off in the solution.}
      \label{fig:disk_origin}
      \end{figure}

To justify the rationale behind half-rank approximations, we tested the algorithm with rank approximations $r=M/8$ and $r=M/1.5$. The results are shown in Figures \ref{Less_than_half_rank} and \ref{More_than_half_rank} respectively. This suggests rank approximations too far away from half-rank approximations can either lead to loss of data or lead to blow-offs which results in improper reconstruction.

\subsubsection{ Test Case 2 - Combined set of phantoms}
\begin{figure}
  \centering
    \subfloat[]{%
     \label{combined_exact} \includegraphics[scale=0.5,keepaspectratio]{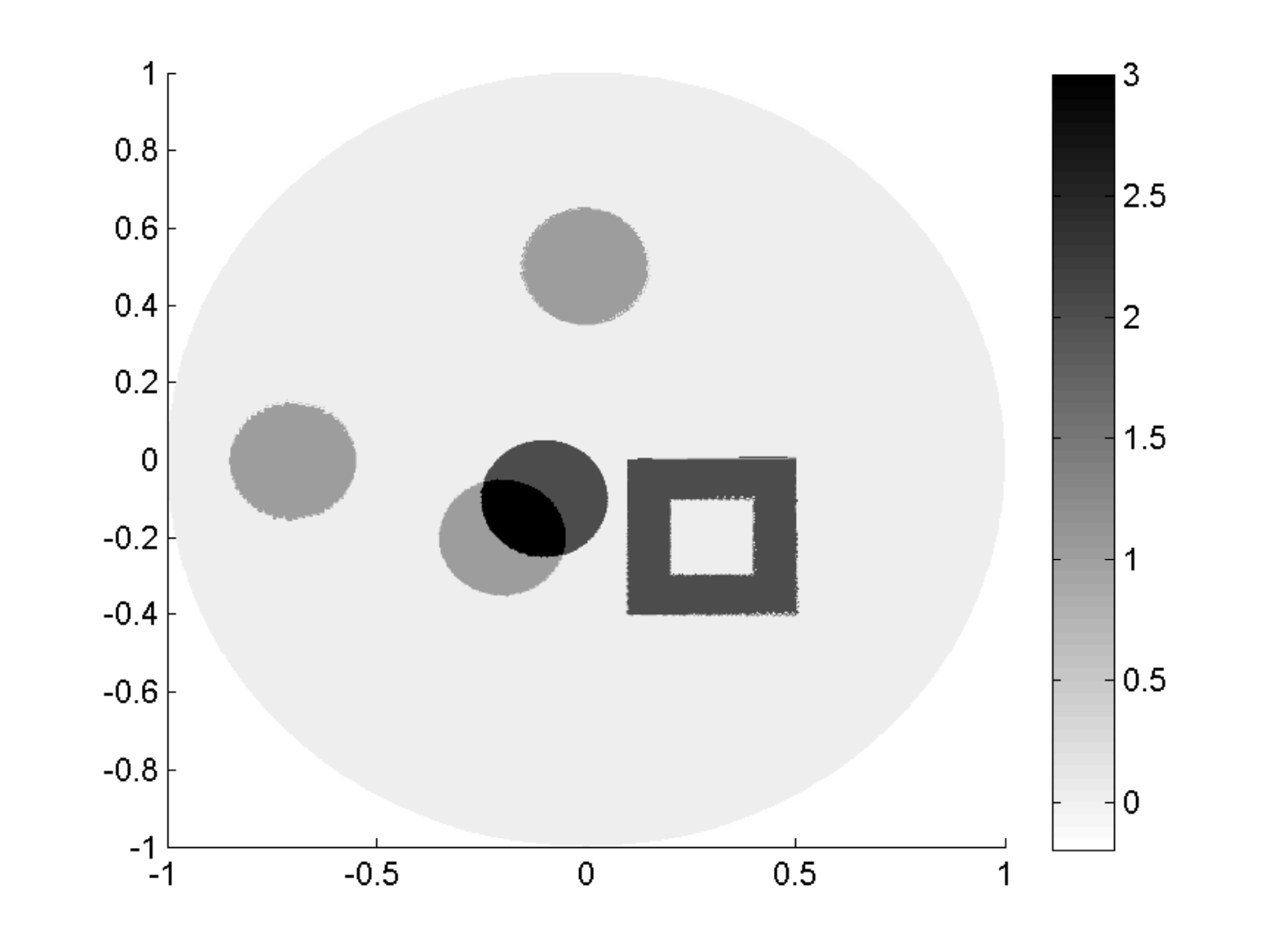}
   }
    \subfloat[]{%
      \label{combined_reconstructed}\includegraphics[scale=0.5,keepaspectratio]{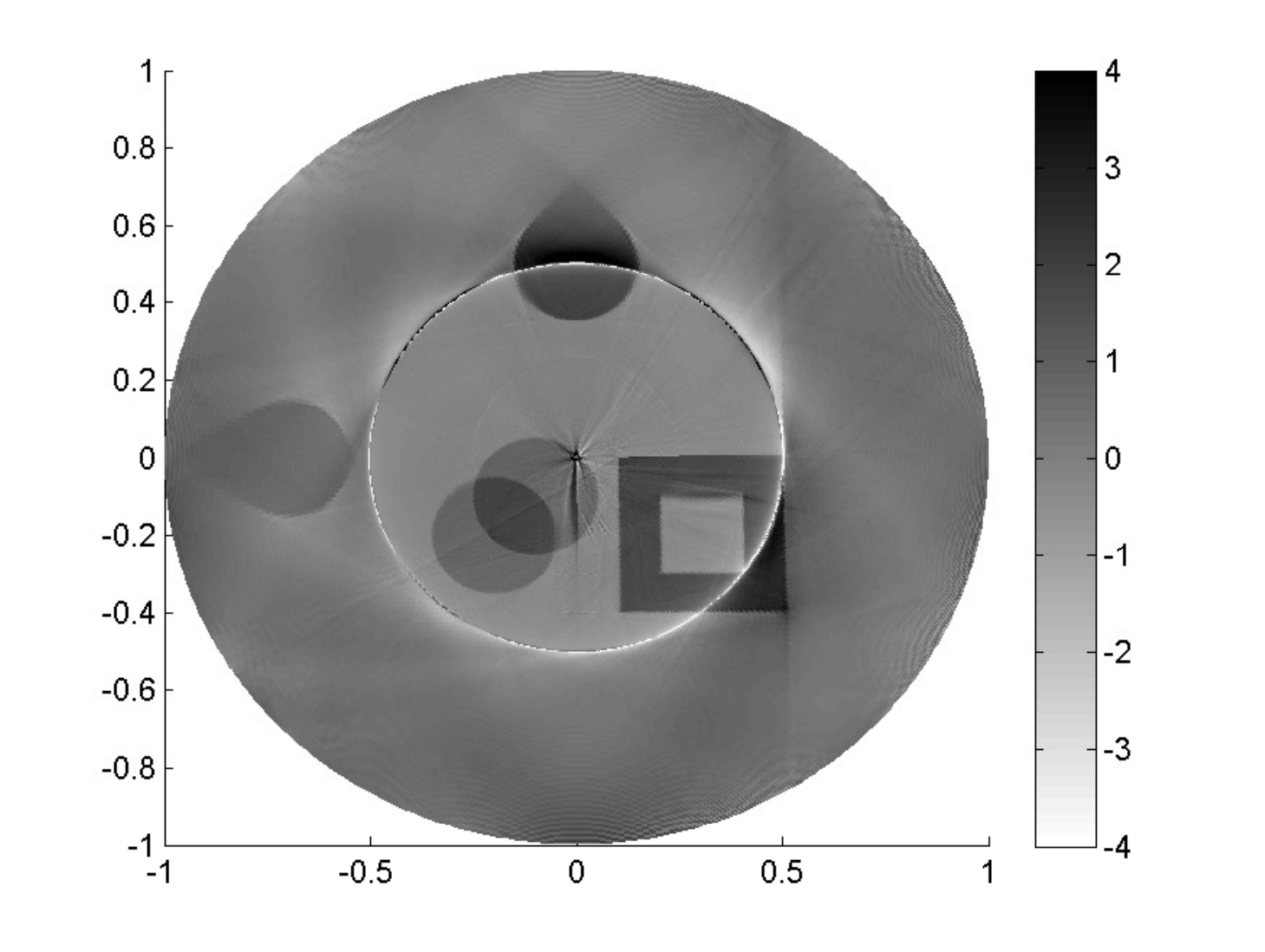}
      }\\
      \subfloat[]{%
      \label{combined_noise}\includegraphics[scale=0.5,keepaspectratio]{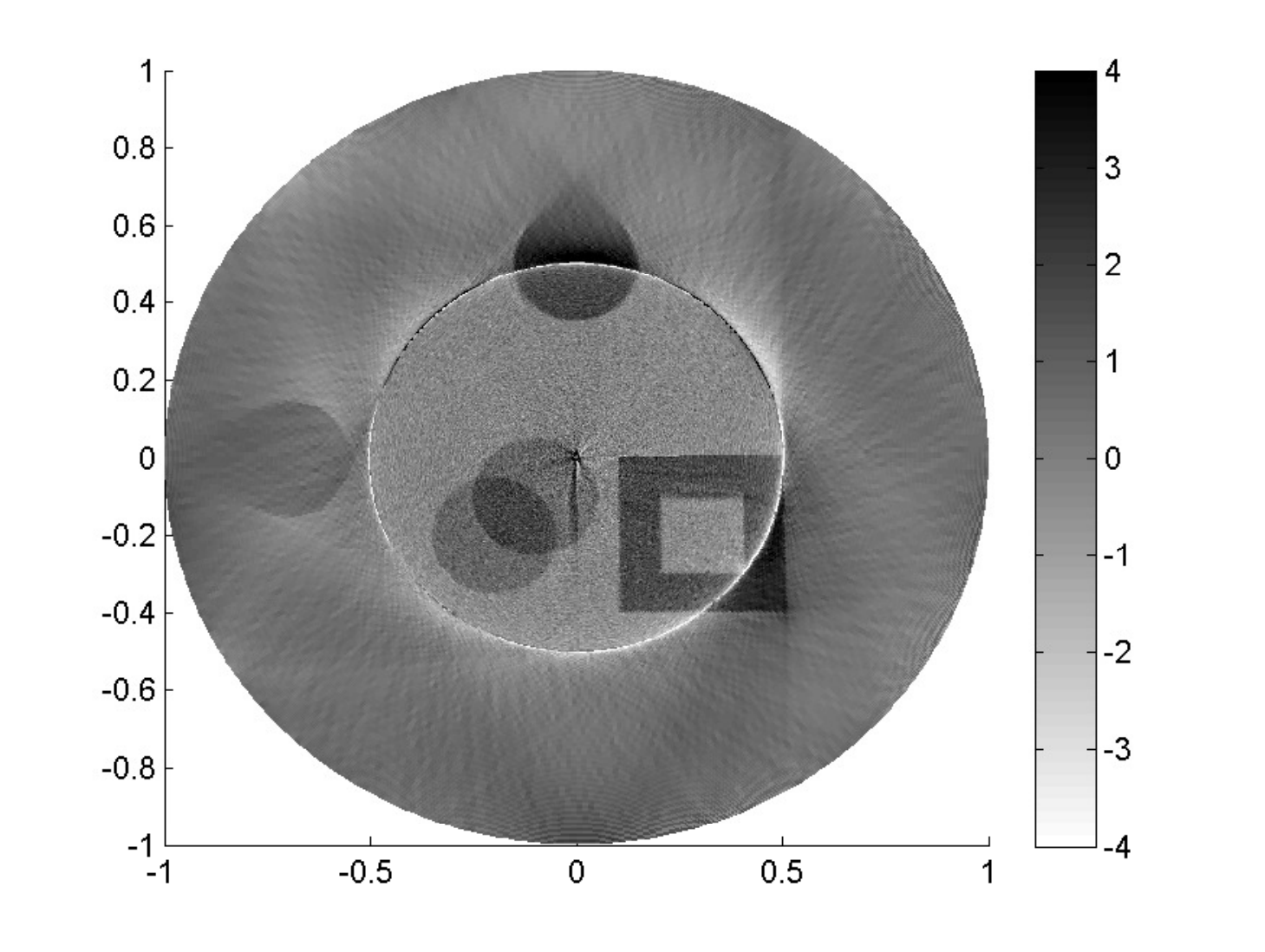}
      }
      \caption{Results for broken ray transform data with breaking angle $\theta=\pi/6$, for a set of phantoms. Figure \ref{combined_exact} represents the actual set of phantoms. Figure \ref{combined_reconstructed} shows the reconstructed image with 400 equally spaced discretizations. Figure \ref{combined_noise} shows the reconstructed image with $5\%$ multiplicative Gaussian noise.}
      \label{fig:combined}
      \end{figure}

In this test case, we consider a set of phantoms represented by a combination of disks with varying intensities and at different locations and a square frame as shown in Fig. \ref{fig:combined}. Fig. \ref{combined_reconstructed} shows the reconstruction with $\theta=\pi/6$ for  400 equally spaced discretizations in $\rho$. Fig. \ref{combined_noise} shows the reconstruction with $5\%$ multiplicative Gaussian noise.

We see in Fig. \ref{combined_reconstructed} and Fig. \ref{combined_noise} that inside the disc of radius $\sin \theta$, we have good reconstructions even with coarser discretizations. Outside the disc of radius $\sin\theta$, we see blurred reconstructions. This is due to lack of stability which will be described next in Sec. \ref{sec:stability_disc}.

\subsection{Stability and Artifacts} \label{sec:stability_disc}
It is a well established fact that image reconstruction in limited data tomography suffers from various types of artifacts. The two most common ones are the blurring of the true singularities of the original object, and the appearance of false singularities that did not exist in the original image (``streak artifacts'') (e.g. see \cite{Frikel-Quinto, Q1993sing}). Here we explain the nature of each of these artifacts, and discuss their appearance in our reconstructed images.

In problems of inverting generalized Radon transforms that integrate a function along smooth curves, one can use standard results of microlocal analysis to predict which parts of the object's (true) singularities will be recovered stably, and which parts will be blurred. In simple words, one can expect to recover stably only those singularities that can be tangentially touched by the integration curves available in the Radon data (e.g. see \cite{XWAK1} and the references there). However, if the generalized Radon transform integrates the image function along non-smooth trajectories, one may be able to do better than that. For example, SSOT in slab geometry produces images of excellent quality using broken rays with basically two (angular) directions (e.g. see \cite{Florescu-Markel-Schotland, Gouia_Amb_V-line}), due to the fact that the integration trajectories can have the scattering ``corner'' at every point of the image domain. In our setup, the reconstructions do not benefit from the presence of corners. The phantom edges do blur if none of the linear pieces of the broken rays touch them tangentially (e.g. see the two discs away from the origin in Fig. \ref{fig:combined} (b) and (c)). The rigorous mathematical study of the ``stabilization due to corners'' or lack of it is not an easy task and is subject of current research by the authors and their collaborators.

The artifacts of second type appear due to the abrupt cut of the (incomplete) Radon data. For example, in limited angle CT the integrals of the image function are available only along lines with limited angular range (e.g. in $[-\pi/4,\pi/4]$ instead of full range $[0,2\pi]$). Then the recovered image may have streak artifacts along lines that have the angular parameters equal to the endpoints of the available limited range. In the CT example above those would be the lines with angular parameters equal to $-\pi/4$ and $\pi/4$ (see \cite{Frikel-Quinto} for more details and a great exposition of this material). In our case, the abrupt cut of the data $\mathcal{R}f(\beta,t)$ happens in two places, when $t=0$ and $t=R$. The first one gives rise to a visible artifact at the origin, since $t=0$ corresponds to rays that break at the origin. The broken rays that correspond to $t=R$ are chords of the disc that pass at distance $R_s=R\sin\theta$ from the origin (see Fig. \ref{Sketch_1}). The envelope of all these chords is the circle of radius $R_s$, along which we have a strong streak artifact in each reconstructed image.

\subsection{Computational Times and Relative $L^2$ error}\label{sec:error}
We now demonstrate the computational efficiency of our developed algorithm by demonstrating the computational times taken and the relative $L^2$ error percentage of reconstruction. The latter is measured only inside the disc of radius $R_s$ with punctured origin to account for errors away from streak artifacts.  All measurements are done for the reconstructions obtained in Sec. \ref{sec:test_cases}.

We define the relative $L^2$ error percentage inside the punctured disc of radius $R_s$ as follows:
$$
\mbox{Relative $L^{2}$ error percentage} = \frac{\|f_{\mbox{rec}}-f_{\mbox{ex}}\|_2}{\|f_{\mbox{ex}}\|_2}*100 \%,
$$
where $f_{\mbox{ex}}=f_{\mbox{ex}}(x_i,y_j)$ and $f_{\mbox{rec}}=f_{\mbox{rec}}(x_i,y_j)$, $i,j=1,\cdots, M$ represents the discretized matrix for the exact function and the reconstructed function respectively, $\|f\|_2=\frac{1}{M^2}\sqrt{\sum_{i=1}^M\sum_{j=1}^Mf_{ij}^2}$ and $f = f(x_i,y_j)$, where $(x_i,y_j)$ lies inside the punctured disc of radius $R_s$.

 The reconstruction of the function $f$ can be divided into two phases
\begin{enumerate}
\item Pre-processing step.
\item Inversion algorithm step.
\end{enumerate}
In the pre-processing step, we compute the inverse of the matrix $A_n$ given in (\ref{mat}) for $n=1,\hdots,N/2+1$ using the half-rank truncated SVD inversion technique described in Section \ref{TSVD}.  The inversion algorithm step consists of the Fourier transform of the Radon data, computing the solution $f_n$ for each $n$, evaluating $f$ using the inverse Fourier transform and finally displaying the results.

\begin{table}[H]
\centering
\begin{tabular}{|c|c|c|c|c|}
\hline
Phantom & Grid&  PP & IA & Error $\%$ \\ [0.5ex]
\hline
Test case 1 & 150 &33.4 sec & 1.1 sec & 35.8\\
\pbox{20mm}{Test case 1\\ (with noise)}&150  & 35.1 sec& 1.2 sec & 36.5 \\
Test case 1 & 400  & 283.4 sec& 7.6 sec & 22.8\\
Test case 1 & 800  & 2988.8 sec& 35.3 sec & 15.8\\
\pbox{20mm}{Test case 2} & 150 & 33.7 sec& 1.4 sec &39.2 \\
[1ex]
\hline
\end{tabular}
\caption{Time taken for the pre-processing step and inversion algorithm for reconstructions with various parameters. PP stands for the pre-processing step and IA stands for the inversion algorithm step. }
\label{time_taken}
\end{table}

Table \ref{time_taken} presents the computational times for the various experiments performed. It can be seen that the computational time depends on the number of radial discretizations of $\rho$ rather than the type of transforms or support of the reconstructed function $f$. We can see from Table \ref{time_taken}, that the computational time taken for the pre-processing step is quite large. But since this step is independent of data, given the number of angular discretizations $N$, this step is computed once, stored in memory and can be used for inversion of any kind of Radon data. This makes the inversion procedure quite fast which can be seen from the time taken for the inversion algorithm step running at less than a minute even for 800 discretizations.

To estimate the computational complexity of our algorithm, we present below the number of operations for various components in the inversion process.
\begin{enumerate}
\item The evaluation of the modified forward FFT for each of the $M$ values of $\rho$ requires $\dfrac{1}{2}~\mathcal{O}(N \log N)$ floating point operations (flops).
\item The next step involves the computation of the truncated SVD for $N/2$ frequencies. For each such frequency, we require $\mathcal{O}(M^3)$ flops.
\item In the next step, we perform the multiplication of the pseudoinverse matrix with the data to compute the Fourier coefficients of the numerical solution. This requires $\mathcal{O}(M^2)$ flops for each of the $N/2$ Fourier coefficients.
\item Finally, we perform the inverse FFT of the obtained Fourier coefficients to get the numerical reconstruction for $M$ values of $\rho$. It requires $\mathcal{O}(N \log N)$ flops.
\end{enumerate}

If the number of discretizations for $\rho$ and $\beta$ are of the same order, i.e. $M$ is of order $N$, then the total number of flops for the inversion process is of order $\mathcal{O}(N^4)$.

\section{Summary}\label{Summary}
We have developed a numerical algorithm for inversion of the broken ray transform in a disk from radially partial data. Our algorithm uses half of the data that the previously known numerical inversions of BRT in the disc used. Given the limitations on the distance that a photon can fly without scattering more than once, our approach allows to double the thickness of objects that can be imaged using single scattering optical tomography. The numerical algorithm requires solution of ill-conditioned matrix problems, which is accomplished using a truncated SVD method. The matrices and the SVD can be constructed in a pre-processing step which can be re-used repeatedly for subsequent computations. This makes our algorithm particularly fast and efficient.
We tested our algorithm on phantoms with jump discontinuities both with and without noise, and it produced high quality reconstructions. The objects in the image were well distinguished, and the recovered intensities of the objects were close to their actual values.

\section{Acknowledgements}\label{Acknowledgements}
The authors thank Venkateswaran P. Krishnan and Eric Todd Quinto for helpful discussions and comments about the paper. 

\bibliographystyle{plain}
\bibliography{references}
\end{document}